\documentclass[11pt]{article}
\renewcommand{\baselinestretch}{1.2}
\usepackage{amsthm,enumerate}
\theoremstyle{plain}
\newtheorem{theorem}{Theorem}[section]
\newtheorem{lemma}[theorem]{Lemma}
\newtheorem{corollary}[theorem]{Corollary}
\newtheorem{proposition}[theorem]{Proposition}
\theoremstyle{definition}
\newtheorem{definition}[theorem]{Definition}
\theoremstyle{remark}
\newtheorem*{remark}{Remark}
\newtheorem*{example}{Example}
\newtheorem*{rems}{Remarks}
\newtheorem*{exs}{Examples}
\newenvironment{remarks}{\begin{rems}\begin{enumerate}}{\end{enumerate}\end{rems}}
\newenvironment{examples}{\begin{exs}\begin{enumerate}}{\end{enumerate}\end{exs}}
\newenvironment{items}{\begin{enumerate}[\rm (i)]}{\end{enumerate}}
\newenvironment{alphitems}{\begin{enumerate}[\rm (a)]}{\end{enumerate}}

 \usepackage{amsmath,amssymb,amsfonts,diagrams}
\newenvironment{keywords}{\noindent\small {\it Keywords\/}:}{\vskip 4pt}
\newenvironment{classification}{\noindent\small 2000 {\it Mathematics Subject
Classification\/}:}{\vskip 12pt}

\newcommand{\comps}{{\mathbb C}}
\newcommand{\reals}{{\mathbb R}}

\newcommand{\posints}{{\mathbb N}}

\newcommand{\free}{{\mathbb F}}

\newcommand{\tensor}{\otimes}

\newcommand{\Tensor}{\hat{\otimes}}

\newcommand{\cstar}{{C^\ast}}

\newcommand{\id}{{\mathrm{id}}}
\newcommand{\cb}{{\mathrm{cb}}}

\newcommand{\A}{{\mathfrak A}}

\newcommand{\Hilbert}{{\mathfrak H}}

\newcommand{\CB}{{\cal CB}}

\newcommand{\VN}{\operatorname{VN}}

\newcommand{\SL}{\operatorname{SL}}

\newcommand{\supp}{{\operatorname{supp}}}

\newcommand{\varcl}[1]{\overline{#1}}

\title{Operator amenability of the Fourier algebra \\
in the $\cb$-multiplier norm}
\author{{\it Brian E.\ Forrest}\thanks{Research supported by NSERC under grant
    no.\ 90749-04.} \and
{\it Volker Runde}\thanks{Research supported by NSERC under grant no.\
  227043-04.} \and 
{\it Nico Spronk}}
\date{}
\begin{document}
\maketitle
\begin{abstract}
Let $G$ be a locally compact group, and let $A_{\cb}(G)$ denote the
 closure of $A(G)$, the Fourier algebra of $G$, in the space of completely 
bounded multipliers of $A(G)$. If $G$ is a weakly  amenable, discrete group 
such that $\cstar(G)$ is residually finite-dimensional, we show that 
$A_{\cb}(G)$ is operator amenable. In particular, 
$A_{\cb}(\free_2)$ is operator amenable even though $\free_2$, the free
group in two generators, is not an amenable group. Moreover, we show that,
if $G$ is a discrete group such that $A_{\cb}(G)$ is operator amenable,
a closed ideal of $A(G)$ is weakly completely complemented in $A(G)$
if and only if it has an approximate identity bounded in the $\cb$-multiplier 
norm.
\end{abstract}
\section*{Introduction}
The Fourier algebra $A(G)$ of a locally compact group $G$ was introduced by P.\ Eymard in \cite{Eym}; for abelian $G$, the Fourier transform
yields an isometric isomorphism of $A(G)$ and $L^1(\hat{G})$, where $\hat{G}$ is the dual group of $G$. Quite soon after the publication of
\cite{Eym}, H.\ Leptin showed that the amenable locally compact groups can be characterized in terms of the Banach algebra $A(G)$ (\cite{Lep}): 
the group $G$ is amenable if and only if $A(G)$ has a bounded approximate identity.
\par
In his memoir \cite{Joh1}, B.\ E.\ Johnson introduced the notion of an amenable Banach algebra and showed that a locally compact group $G$
is amenable if and only if $L^1(G)$ is amenable (this result motivates the choice of terminology). Since every amenable Banach algebra has
a bounded approximate identity, Leptin's theorem immediately yields that, for any locally compact group $G$, the amenability of $A(G)$ 
necessitates that of $G$. The question for which locally compact groups $G$ precisely the Fourier algebra $A(G)$ is amenable remained open
for a surprisingly long period of time. In \cite{Joh3}, Johnson showed that $A(G)$ may fail to be amenable for certain compact groups $G$.
Eventually, the first- and the second-named author proved that $A(G)$ is amenable if and only if $G$ has an abelian subgroup of finite index
(\cite{FR}).
\par
The study of the Fourier algebra gained new momentum in 1995 with \cite{Rua}. As the predual of the group von Neumann algebra,
$A(G)$ is an operator space in a canonical manner for every locally compact group $G$. Z.-J.\ Ruan used this to add operator space overtones
to Johnson's notion of an amenable Banach algebra and introduced the concept of operator amenability. In \cite{Rua}, he showed that a locally
compact group $G$ is amenable if and only if $A(G)$ is operator amenable. Since then, it has become apparent that the theory of operator spaces
(\cite{ER}) provides powerful tools for the study of the Fourier algebra of a
locally compact group and of related algebras (\cite{Ari}, \cite{ARS},
\cite{FW}, \cite{IS}, \cite{LNR}, \cite{RunPP}, \cite{RS}, \cite{Spr}, \cite{Woo1}, and \cite{Woo2}; for a detailed overview, see \cite{RunExpo}). 
Even if one is mainly interested in the Fourier algebra as a mere Banach algebra, operator space methods provide new insights: the main 
results of both \cite{FKLS} and \cite{FR} do not make any reference to operator spaces, but their respective proofs depend on operator space 
techniques.
\par
In this paper, we investigate another Banach algebra associated with a 
locally compact group $G$, mostly from an operator space point of view.
The Fourier algebra $A(G)$ embeds canonically into the algebra of 
completely bounded multipliers on $A(G)$. For amenable $G$, the norm 
on $A(G)$ inherited from this algebra --- the $\cb$-multiplier norm --- 
is the given norm; for non-amenable groups, however, the two norms are
inequivalent. We denote the completion of $A(G)$ with respect to the
$\cb$-multiplier norm by $A_{\cb}(G)$. Unlike the Fourier algebra, $A_{\cb}(G)$
may well have a bounded approximate identity for non-amenable $G$ (such 
groups are called weakly amenable). The main result of this paper is that, 
for certain discrete, non-amenable groups --- among them $\free_2$, the
free groups in two generators ---, $A_{\cb}(G)$ not only has a bounded
approximate identity, but is operator amenable. We then move on to study
complementation properties of closed ideals of both $A_{\cb}(G)$ and $A(G)$, 
where $G$ is discrete with $A_{\cb}(G)$ operator amenable. In particular, we 
show that, for such $G$, a closed ideal of $A(G)$ is (weakly) completely 
complemented in $A(G)$ if and only if it has an approximate identity that
is bounded in $A_{\cb}(G)$.
\section{Preliminaries}
Our reference for the theory of operator spaces is the monograph \cite{ER}, whose notation and choice of terminology we adopt unless explicitly
stated otherwise.
\par
We begin with introducing basic definitions:
\begin{definition}
A \emph{quantized Banach algebra} is an algebra $\A$ which is also an operator space such that the multiplication of $\A$ is completely 
bounded.
\end{definition}
\begin{remark}
We do not require the multiplication of a quantized Banach to be completely contractive: this extra bit of generality can be convenient 
sometimes (\cite{LNR}).
\end{remark}
\begin{examples}
\item If $\A$ is any Banach algebra, then $\max \A$ (the maximal operator space over $\A$; see \cite{ER}) is a quantized Banach algebra.
\item If $\Hilbert$ is a Hilbert space, then every closed subalgebra of 
${\mathcal B}(\Hilbert)$ with its concrete operator space structure is
a quantized Banach algebra.
\item Let $E$ be an operator space. Then $\CB(E)$ is a quantized Banach 
algebra.
\item Let $G$ be a locally compact group, let $\VN(G)$ denote its group von Neumann algebra, and let $\cstar(G)$ and $C^\ast_r(G)$ denote
its full and reduced group $\cstar$-algebra, respectively. The dualities
\[
  A(G) = \VN(G)_\ast, \qquad B(G) = \cstar(G)^\ast, \qquad\textnormal{and}\qquad B_r(G) = C^\ast_r(G)^\ast 
\]
equip $A(G)$ as well as $B(G)$ and $B_r(G)$, the {\it Fourier--Stieltjes algebra\/} and the {\it reduced Fourier--Stieltjes algebra\/} 
(\cite{Eym}), respectively, with an operator space structure. With these operator space structures, $A(G)$, $B(G)$, and $B_r(G)$ are
quantized Banach algebras.
\item Let $G$ be a locally compact group. A \emph{multiplier} of $A(G)$ is 
a (necessarily bounded and continuous) function $f \!: G \to \mathbb{C}$ such 
that $f A(G) \subset A(G)$. For each such $f$, multiplication with $f$ is a 
linear operator on $A(G)$ --- bounded by the closed graph theorem --- which we 
denote by $M_f$; it is straightforward that $M_f \!: A(G) \to A(G)$ is an 
$A(G)$-module homomorphism. Alternatively, the term multiplier is also
used to refer to an $A(G)$-module homomorphism on $A(G)$. Both usages are
equivalent: whenever $T \!: A(G) \to A(G)$ is an $A(G)$-module homomorphism,
there is (a necessarily unique) $f \!: G \to \comps$ with $f A(G) \subset 
A(G)$ such that $T = M_f$ (\cite[p.\ 422]{Dal}). The 
\emph{multiplier algebra} of $A(G)$ is the closed subalgebra
\[
  {\mathcal M}(A(G)) := \{ M_f : \text{$f$ is a multiplier of $A(G)$} \}
\]
of ${\mathcal B}(A(G))$. For notational convenience, we shall simply 
identify a multiplier $f$ of $A(G)$ and the corresponding $M_f$.
Finally, the \emph{$\cb$-multiplier algebra} of $A(G)$ is 
defined as
\[
  {\mathcal M}_{\cb}(A(G)) := \CB(A(G)) \cap {\mathcal M}(A(G));
\]
it is a closed subalgebra of $\CB(A(G))$ and thus a quantized Banach
algebra. 
\end{examples}
\begin{definition} \label{qmod}
Let $\A$ be a quantized Banach algebra. A {\it quantized Banach $\A$-bimodule\/} is an $\A$-bimodule $E$ which is also an operator space such 
that the module actions
\[
  \A \times E \to E, \quad (a,x) \rightarrow a \cdot x
  \qquad\textnormal{and}\qquad
  E \times \A \to E, \quad (x,a) \rightarrow x \cdot a
\]
are completely bounded.
\end{definition}
\begin{remark}
Our quantized Banach bimodules are not to be confused with the operator 
bimodules studied, for instance, in \cite{BL}: every operator bimodule over an
operator algebra is a quantized Banach bimodule in the sense of Definition
\ref{qmod}, but the converse is false.
\end{remark}
\par
If $\A$ is a quantized Banach algebra and $E$ is a quantized Banach 
$\A$-bimodule, then $E^\ast$ becomes a quantized Banach $\A$-bimodule 
in a canonical way through
\[
  \langle x, a \cdot \phi \rangle := \langle x \cdot a, \phi \rangle
  \quad\textnormal{and}\quad
  \langle x, \phi \cdot a \rangle := \langle x, a \cdot \phi \rangle
  \qquad (a \in \A, \, \phi \in E^\ast, \, x \in E).
\]
\begin{definition} \label{opamdef}
A quantized Banach algebra $\A$ is said to be \emph{operator amenable} if, for every quantized Banach $A$-bimodule $E$, every completely
bounded derivation $D \!: \A \to E^\ast$ is inner.
\end{definition}
\begin{examples}
\item Let $G$ be a locally compact group. Then $A(G)$ is operator amenable if and only if $G$ is amenable (\cite[Theorem 3.6]{Rua}).
\item Let $\A$ be a Banach algebra. Then the quantized Banach algebra $\max \A$ is amenable if and only if $\A$ is amenable (\cite{Rua})
\item A $\cstar$-algebra is amenable if and only if it is operator amenable (\cite[Theorem 5.1]{Rua}).
\end{examples}
\par
We also require a modification of Definition \ref{opamdef}.
\par
In \cite{RunStudia}, the second-named author considered a class of Banach algebras --- suggestively named \emph{dual Banach algebras} ---
which are dual Banach space (with a fixed, but not necessarily unique predual) such that multiplication is separately $w^\ast$-continuous.
In \cite{RS}, the second- and the third-named author extended this notion to a quantized setting:
\begin{definition}
A quantized Banach algebra $\A$ is called \emph{dual} if $\A = (\A_\ast)^\ast$ for some Banach space $\A_\ast$ such that the multiplication
of $\A$ is separately $\sigma(\A,\A_\ast)$-continuous.
\end{definition}
\begin{examples}
\item If $\A$ is a dual Banach algebra in the sense of \cite{RunStudia}, then $\max \A$ is a dual, quantized Banach algebra.
\item Every von Neumann algebra is a dual, quantized Banach algebra.
\item Let $G$ be a locally compact group. Then $B(G)$ and $B_r(G)$ are dual, quantized Banach algebras (\cite{RS}).
\item Let $G$ be a locally compact group. In both \cite{dCH} and \cite{Spr2}, a predual space of ${\mathcal M}_{\cb}(A(G))$ is constructed.
A priori, it is not clear that these two predual spaces are identical. However, from \cite[Lemma 1.9]{dCH} and \cite[Corollary 6.6]{Spr2},
it follows that, on norm bounded subsets of ${\mathcal M}_{\cb}(A(G))$, the $w^\ast$-topology on ${\mathcal M}_{\cb}(A(G))$ induced by either predual space 
is the relative topology of $\sigma(L^\infty(G),L^1(G))$. The Krein--\v{S}mulian theorem then yields that the predual spaces from
both \cite{dCH} and \cite{Spr2} are identical. Since multiplication in $L^\infty(G)$ is separately $\sigma(L^\infty(G),L^1(G))$-continuous,
we obtain that multiplication in ${\mathcal M}_{\cb}(A(G))$ is separately $w^\ast$-continuous, first on norm-bounded sets, and then --- by virtue
of the Krein--\v{S}mulian theorem again --- on all of ${\mathcal M}_{\cb}(A(G))$. Hence, ${\mathcal M}_{\cb}(A(G))$ is a dual, quantized Banach algebra.
\end{examples}
\par
In \cite{RunStudia}, a weaker variant of amenability --- dubbed Connes-amenability --- was introduced for dual Banach algebras. Generally,
Connes-amenability seems to be better suited for dual Banach algebras than the
original definition from \cite{Joh1} (compare \cite{DGH} and \cite{RunLMS}, for
example). In \cite{RS}, the second- and the third-named author extended the notion of Connes-amenability to the quantized setting.
\par
Let $\A$ be a quantized Banach algebra, and let $E$ be a dual, quantized Banach $\A$-bimodule, i.e.\ a quantized Banach $\A$-bimodule which
is the canonical dual module of some other quantized Banach $\A$-bimodule. Suppose that $\A$ is dual. Then we say that $E$ is \emph{normal} if
the module actions
\[
  \A \times E \to E, \quad (a,x) \rightarrow a \cdot x
  \qquad\textnormal{and}\qquad
  E \times \A \to E, \quad (x,a) \rightarrow x \cdot a
\]
are separately $w^\ast$-continuous.
\begin{definition}
A dual, quantized Banach algebra $\A$ is said to be {\it operator Connes-amenable\/} if, for every normal, dual, quantized Banach 
$\A$-bimodule $E$, every $w^\ast$-con\-ti\-nu\-ous, completely bounded derivation $D \!: \A \to E$ is inner.
\end{definition}
\begin{examples}
\item A dual Banach algebra $\A$ is Connes-amenable in the sense of \cite{RunStudia} if and only if $\max \A$ is operator Connes-amenable.
\item A locally compact group $G$ is amenable if and only if $B_r(G)$ is operator Connes-amenable (\cite[Theorem 4.4]{RS}).
\item The free group in two generators, which we denote by $\free_2$, is not amenable, but $B(\free_2)$ is operator Connes-amenable
(\cite{RS}).
\end{examples}
\section{Operator amenability of ${A_{\cb}}(G)$ for non-amenable $G$}
Let $G$ be a locally compact group. Then we have the following completely contractive inclusions:
\[
  A(G) \subset B_r(G) \subset B(G) \subset {\mathcal M}_{\cb}(G).
\]
The first and the second inclusion are always complete isometries, whereas the
embedding of $A(G)$ into ${\cal M}_{\cb}(A(G))$ is bounded below only
if $G$ is amenable. In this case, we have completely isometric 
identifications $B_r(G) = B(G) = {\cal M}_{\cb}(A(G))$, so that $A(G)$ 
embeds into ${\cal M}_{\cb}(G)$ completely isometrically. For a discussion 
and further references, see \cite{Spr2}.
\par
As ${\mathcal M}_{\cb}(A(G))$ is a dual, quantized Banach algebra, it makes sense to ask for which locally compact groups $G$ it is operator
Connes-amenable. Of course, if $G$ is amenable, then ${\cal M}_{\cb}(A(G))
= B(G) = B_r(G)$ is operator Connes-amenable (\cite[Theorem 4.4]{RS}).
The following proposition, gives another sufficient condition.
\par
Recall that a $\cstar$-algebra is said to be 
\emph{residually finite-dimensional} if its fi\-nite-di\-men\-sio\-nal, 
irreducible $^\ast$-representation separate its points. Furthermore, 
following \cite{HK}, we say that a locally compact group $G$ has the 
\emph{approximation property} if there is a net in $A(G)$ converging to
the constant function $1$ in the $w^\ast$-topology of ${\cal M}_{\cb}(A(G))$.
\begin{proposition} \label{prop1}
Let $G$ be a locally compact group with the approximation property
such that $\cstar(G)$ is residually finite-dimensional. 
Then ${\mathcal M}_{\cb}(A(G))$ is operator Connes-am\-en\-ab\-le.
\end{proposition}
\begin{proof}
Since $G$ has the approximation property and since ${\cal M}_{\cb}(A(G))$ is
a dual Banach algebra, it is clear that $A(G)$ is 
$w^\ast$-dense in ${\cal M}_{\cb}(A(G))$. Consequently, $B(G) \supset A(G)$ is 
also $w^\ast$-dense in ${\cal M}_{\cb}(A(G))$.
Since $\cstar(G)$ is residually finite-dimensional, $B(G)$ 
is operator Connes-amenable by \cite[Theorem 4.6]{RS}. As remarked earlier,
the $w^\ast$-topologies on both $B(G)$ and ${\cal M}_{\cb}(A(G))$ coincide 
on bounded sets with the relative topology induced by 
$\sigma(L^\infty(G),L^1(G))$, so that the inclusion $B(G) \subset
{\cal M}_{\cb}(A(G))$ is $w^\ast$-continuous by the Krein--\v{S}mulian theorem.
From (the quantized analog of) \cite[Proposition 4.2]{RunStudia}, it then
follows that ${\cal M}_{\cb}(A(G))$ is also operator Connes-amenable.
\end{proof}
\begin{example}
By \cite[Proposition VII.6.1]{Dav}, $\cstar(\free_2)$ is residually 
finite-dimensional, and, as we shall note below, $\free_2$ has the 
approximation property. Hence, ${\mathcal M}_{\cb}(A(\free_2))$
is operator Connes-amenable by Proposition \ref{prop1} --- even though 
$\free_2$ fails to be amenable.
\end{example}
\par
If $G$ is an amenable, locally compact group, then $A(G)$ embeds completely 
isometrically into ${\mathcal M}_{\cb}(A(G)) = B(G)$. If $G$ is 
not amenable, however, $A(G)$ is not closed in ${\mathcal M}_{\cb}(A(G))$.
We convene to denote the closure of $A(G)$ in ${\mathcal M}_{\cb}(A(G))$ by 
$A_{\cb}(G)$. 
\par
We collect a few basic properties of $A_{\cb}(G)$:
\begin{proposition}
Let $G$ be a locally compact group. Then ${A_{\cb}}(G)$ is a regular,
commutative, Tauberian (quantized)
Banach algebra whose character space is canonically 
identified with $G$.
\end{proposition}
\begin{proof}
It is clear that $A_{\cb}(G)$ is commutative and semisimple, and
by \cite[Lemma 1]{For1}, the character space of $A_{\cb}(G)$ is $G$
in the canonical way. Let $F \subset G$ be closed, and let $x \in G
\setminus F$. Since $A(G)$ is regular (see \cite{Eym}), there is $f \in A(G)
\subset A_{\cb}(G)$ such that $f |_F \equiv 0$ and $f(x) = 1$.
Hence, $A_{\cb}(G)$ is also regular.
\par
To see that $A_{\cb}(G)$ is Tauberian, let $f \in A_{\cb}(G)$,
and let $\epsilon > 0$. Since $A(G)$ is dense in $A_{\cb}(G)$, there is
$g \in A(G)$ with $\| f - g \|_{A_{\cb}(G)} < \frac{\epsilon}{2}$, and
since $A(G)$ is Tauberian, there is $h \in A(G)$ with $\supp(h)$ compact
and
\[
  \| g - h \|_{A_{\cb}(G)} \leq \| g - h \|_{A(G)} < \frac{\epsilon}{2}.
\]
It follows that $\| f - h \|_{A_{\cb}(A)} < \epsilon$. Since $\epsilon > 0$ is arbitrary, 
this means that $A_{\cb}(G)$ is Tauberian.
\end{proof}
\par
Leptin's theorem (\cite{Lep}), motivates the adverb ``weakly'' in the following definition:
\begin{definition} \label{wam}
A locally compact group is said to be \emph{weakly amenable} if ${A_{\cb}}(G)$ has a bounded approximate identity.
\end{definition}
\begin{remarks}
\item This definition of a weakly amenable, locally compact group isn't quite the original one (compare \cite{CH}), but is easily seen to be equivalent 
(\cite[Proposition 1]{For}).
\item Weakly amenable groups have the approximation property of \cite{HK}
whereas the converse is false.
\item In \cite{BCD}, a notion of weak amenability for Banach algebras was
introduced. This Banach algebraic amenability --- which can easily be
adapted to the quantized setting --- is related to Definition \ref{wam} only
in the sense that both weak amenabilities are weaker than the notions of
amenability for Banach algebras and locally compact groups, respectively. 
There are no analogs of \cite[Theorem 2.5]{Joh1} or \cite[Theorem 3.6]{Rua}: 
$L^1(G)$ is weakly amenable (\cite{Joh2}) and $A(G)$ is operator weakly 
amenable (\cite{Spr}) for every locally compact group $G$.
\end{remarks}
\begin{examples}
\item By \cite{Lep}, every amenable, locally compact groups is weakly amenable.
\item Even though $\free_2$ is not amenable, it is weakly amenable (\cite[Corollary 3.9]{dCH}).
\item It is shown in \cite{Dor1} that $\SL(2,\reals) \rtimes \reals^N$ is not
weakly amenable for $N \geq 2$. In \cite{Dor2}, this is used to show that
every simple Lie group with real rank greater than or equal to two fails to
be weakly amenable.
\end{examples}
\par
It is clear from \cite[Theorem 3.6]{Rua} --- combined with elementary hereditary properties of operator amenability --- that
${A_{\cb}}(G)$ is operator amenable for every amenable, locally compact group. In the remainder of this section, we shall see that the converse 
need not be true. 
\par
We first present three lemmas.
\par
Let $\A$ be a Banach algebra, and recall that a Banach $\A$-bimodule is called \emph{pseudo-unital} (or \emph{neo-unital}) if
\[
  E = \{ a \cdot x \cdot b : a,b \in \A, \, x \in E \}.
\]
\begin{lemma} \label{lem1}
Let $\A$ be a quantized Banach algebra with a bounded approximate identity. 
Then $\A$ is operator amenable if and only if, for
each pseudo-unital, quantized Banach $\A$-bimodule $E$, every completely 
bounded derivation $D \!: \A \to E^\ast$ is inner.
\end{lemma}
\begin{proof} 
The proof for the classical case (\cite[Proposition 2.1.5]{LoA}) carries over nearly verbatim.
\end{proof}
\par
Let $\A$ be a Banach algebra, and let $I$ be a closed ideal of $\A$. The 
\emph{$I$-strict} topology on $\A$ is the locally convex topology
induced by the seminorms
\[
  \A \to [0,\infty), \quad a \mapsto \| ax \| + \| xa \| \qquad (x \in I).
\]
(Note that this topology need not be Hausdorff.)
\begin{lemma} \label{lem2}
Let $\A$ be a quantized Banach algebra, let $I$ be a closed ideal of $\A$ 
with a bounded approximate identity, let $E$ be a pseudo-unital,
quantized Banach $I$-bimodule, and let $D \!: I \to E^\ast$ be a completely
 bounded derivation. Then $E$ is a quantized Banach $\A$-bimodule
in a canonical fashion, and there is a completely bounded derivation 
$\tilde{D} \!: \A \to E^\ast$ extending $D$ which is continuous with
respect to the $I$-strict topology on $\A$ and the $w^\ast$-topology on 
$E^\ast$.
\end{lemma}
\begin{proof}
By \cite[Proposition 2.1.6]{LoA}, the module action of $I$ on $E$ extends
canonically to $\A$, and $D$ has a bounded extension $\tilde{D} \!: \A 
\to E^\ast$ which is continuous with respect to the $I$-strict topology on 
$\A$ and the $w^\ast$-topology on $E^\ast$. (Since $I$ is dense in $\A$ in
the $I$-strict topology, $\tilde{D}$ is uniquely determined by its 
continuity properties.)
\par
Two claims remain to be checked:
that $E$ is indeed a \emph{quantized} Banach $\A$-bimodule, and that
$\tilde{D}$ is \emph{completely} bounded.
\par
We first verify that $E$ is a quantized Banach $\A$-bimodule.
Let $( e_\alpha )_\alpha$ be an approximate identity for $I$ bounded by
$C \geq 0$. Note that, since $E$ is a pseudo-unital Banach $\A$-bimodule, we 
have
\[
  \lim_\alpha [ e_\alpha \cdot x_{j,k} ] = [ x_{j,k} ] =
  \lim_\alpha [ x_{j,k} \cdot e_\alpha ] \qquad (n \in \posints, \, 
  [x_{j,k}] \in M_n(E)).
\]
Let $\kappa \geq 0$ be the $\cb$-norm of the completely bounded, bilinear map
$I \times E \ni (b,x) \mapsto b \cdot x$, and fix $n \in \posints$. Let
$[a_{j,k}] \in M_n(\A)$ and let $[ x_{\nu,\mu} ] \in M_n(E)$. 
It follows that
\[
  \begin{split}
  \| [a_{j,k}\cdot \xi _{\nu,\mu} ] \|_{M_{n^2}(E)} & = 
  \lim_\alpha 
  \|[a_{j,k} \cdot (e_\alpha \cdot x _{\nu,\mu})]\|_{M_{n^2}(E)} \\ 
  & = \lim_\alpha \|[a_{j,k} e_\alpha \cdot x _{\nu,\mu}]\|_{M_{n^2}(E)} \\
  & \leq \kappa \limsup_\alpha 
  \|[a_{j,k} e_\alpha ]\|_{M_n(I)} \| [x_{\nu,\mu}]\|_{M_n(E)} \\
  & \leq \kappa C \|[a_{j,k}]\|_{M_n(\A)} \| \| [x_{\nu,\mu}]\|_{M_n(E)} ,\\
\end{split}
\]
so that the extended module left module action 
action $\A \times E \ni (a,x) \mapsto a \cdot x$ is completely bounded (by
$\kappa C$). Similarly, one sees that $E \times \A \ni (x,a) \mapsto x \cdot
a$ is completely bounded. Consequently, $E$ is indeed a quantized Banach
$\A$-bimodule.
\par
Next, we turn to showing that the extension $\tilde{D} \!: \A \to E^\ast$
from \cite[Proposition 2.1.6]{LoA} is not only bounded, but completely bounded.
Let $a \in \A$, let $x \in E$, and let $b,c \in I$, and note that
\begin{multline*}
  \left\langle b \cdot x \cdot c, \tilde{D}a \right\rangle \\
  =
  \lim_\alpha 
  \left\langle e_\alpha b \cdot x \cdot c, \tilde{D} a \right \rangle =
  \lim_\alpha 
  \left\langle b \cdot x \cdot c, \left( \tilde{D} a \right) \cdot e_\alpha
  \right \rangle =
  \lim_\alpha 
  \langle b \cdot x \cdot c, 
  D (a e_\alpha) - a \cdot D(e_\alpha) \rangle.
\end{multline*}
Since $E$ is pseudo-unital, this means that
\[
  \tilde{D} a = \text{$\sigma(E^\ast,E)$-}\lim_\alpha 
  (D (a e_\alpha) - a \cdot D(e_\alpha))
  \qquad (a \in \A)
\]
and, consequently,
\begin{multline*}
  \tilde{D}^{(n)}([ a_{j,k}]) \\ = 
  \text{$\sigma(M_n(E^\ast),T_n(E))$-}\lim_\alpha 
  (D ([a_{j,k} e_\alpha]) - [a_{j,k} \cdot D(e_\alpha)])
  \qquad (n \in \posints, \, [a_{j,k}] \in M_n(\A)),
\end{multline*}
where, $\tilde{D}^{(n)} \!: M_n(\A) \to M_n(E^\ast)$ denotes the $n$-th
amplification of $\tilde{D}$ for $n \in \posints$. To see that $\tilde{D}$
is completely bounded, let $n \in \posints$ and $[a_{j,k}]\in M_{n}(A)$,
and note that, by the foregoing,
\[
  \begin{split}
  \left\| \tilde{D}^{(n)}([a_{j,k}]) \right\|_{M_n(E^\ast)} & \leq
  \limsup_\alpha 
  \| (D ([a_{j,k} e_\alpha]) - [a_{j,k} \cdot D(e_\alpha)]\|_{M_n(E^\ast)} \\
  & \leq \limsup_\alpha ( \| D \|_{\cb} \| a_{j,k} \|_{M_n(\A)} \| e_\alpha \|+
  \tilde{\kappa} \| a_{j,k} \|_{M_n(\A)} \| D \| \| e_\alpha \| ) \\
  & \leq ( C \| D \|_{\cb} + \tilde{\kappa} C \| D \|) \| a_{j,k} \|_{M_n(\A)},
  \end{split}
\]
where $\tilde{\kappa}$ is the $\cb$-norm of the left module action 
$\A \times E^\ast \ni (a,\phi) \mapsto a \cdot \phi$.
Hence, $\tilde{D}$ is indeed completely bounded (with $\left\| \tilde{D}
\right\|_{\cb} \leq C \| D \|_{\cb} + \tilde{\kappa} C \| D \|$).
\end{proof}
\par
Our final lemma is:
\begin{lemma} \label{lem3}
Let $G$ be a discrete group. Then the following topologies coincide on norm bounded subsets of ${\mathcal M}_{\cb}(A(G))$:
\begin{alphitems}
\item the $w^\ast$-topology;
\item the topology of pointwise convergence on $G$;
\item the ${A_{\cb}}(G)$-strict topology.
\end{alphitems}
\end{lemma}
\begin{proof}
That (a) and (b) coincide on norm bounded subsets follows from 
\cite[Lemma 1.9]{dCH}, and the fact that ${A_{\cb}}(G)$ is Tauberian 
yields the corresponding statement for (b) and (c).
\end{proof}
\par
We can now state and prove the main result of this section:
\begin{theorem} \label{thm1}
Let $G$ be a weakly amenable, discrete group such that $\cstar(G)$ is 
residually finite-dimensional. Then ${A_{\cb}}(G)$ is operator amenable.
\end{theorem}
\begin{proof}
Let $E$ be a quantized Banach $A_{\cb}(G)$-bimodule, and let 
$D \!: A_{\cb}(G) \to E^\ast$ be a completely bounded derivation.
Since $G$ is weakly amenable, i.e.\ $A_{\cb}(G)$ has a bounded 
approximate identity, we may invoke Lemma \ref{lem1} and suppose without 
loss of generality that $E$ is pseudo-unital. By Lemma \ref{lem2}, $E$ is a 
quantized Banach ${\mathcal M}_{\cb}(A(G))$-bimodule in a canonical
way, and there is a completely bounded derivation $\tilde{D} \!: 
{\mathcal M}_{\cb}(A(G)) \to E^\ast$ that extends $D$ and is continuous with
respect to the $A_{\cb}(G)$-strict topology on 
${\mathcal M}_{\cb}(A(G))$ and the $w^\ast$-topology on $E^\ast$. 
\par
Due to Lemma \ref{lem3}, an argument as in the proof of 
\cite[Theorem 3.5]{RunScand} yields that the dual, quantized Banach 
${\mathcal M}_{\cb}(A(G))$-module $E^\ast$ is actually normal and that 
$\tilde{D}$ is $w^\ast$-$w^\ast$-continuous.
\par
Since $\cstar(G)$ is residually finite-dimensional and since $G$
has the approximation property, ${\mathcal M}_{\cb}(A(G))$ 
is operator Connes-amenable by Proposition \ref{prop1}.
Consequently, $\tilde{D}$ --- and therefore $D$ --- is inner.
\end{proof}
\par
With Theorem \ref{thm1} proven, it is not hard to come up with examples of locally compact groups $G$ that fail to be amenable, but for
which ${A_{\cb}}(G)$ is nevertheless operator amenable:
\begin{example}
Since $\free_2$ is weakly amenable and $\cstar(\free_2)$ is residually
finite-dimensional, ${A_{\cb}}(\free_2)$ is operator amenable 
by Theorem \ref{thm1}.
\end{example}
\par
Even though we have exhibited non-amenable (discrete) groups $G$ 
for which $A_{\cb}(G)$ is operator amenable, we are still far from a 
characterization of those locally compact groups $G$ such that $A_{\cb}(G)$ 
is operator amenable. It may be that $A_{\cb}(G)$
is operator amenable whenever $G$ is weakly amenable.
\par
As in \cite{ER}, $\Tensor$ stands for the projective tensor product of 
operator spaces. If $\A$ is a quantized Banach algebra, $\A \Tensor \A$
becomes a quantized Banach $\A$-bimodule via
\[
  a \cdot (x \tensor y) := ax \tensor y \quad\text{and}\quad
  (x \tensor y) \cdot a := x \tensor ya \qquad (a,x,y \in \A),
\]
so that the multiplication operator
\[
  \Delta \!: \A \Tensor \A \to \A, \quad a \tensor b \mapsto ab
\]
becomes a completely bounded homomorphism of $\A$-bimodules.
\par
The following definition arises naturally in A.\ Ya.\ Helemski\u{\i}'s
topological homology (\cite{Hel}) --- or rather in its quantized version
(see \cite{Ari} or \cite{Hel2}, for example):
\begin{definition}
A quantized Banach algebra $\A$ is called \emph{operator 
biprojective} if the multiplication operator 
$\Delta \!: \A \Tensor \A \to \A$ has a completely
bounded right inverse which is also a homomorphism of $\A$-bimodules.
\end{definition}
\begin{example}
Let $G$ be a locally compact group. As was shown independently by O.\ Yu.\
Aristov (\cite{Ari}) and P.\ J.\ Wood (\cite{Woo2}), $A(G)$ is operator 
amenable if and only if $G$ is discrete.
\end{example}
\par
Let $G$ be any locally compact group such that $A_{\cb}(G)$ is
operator biprojective. Then $G$ has to be discrete by (the quantized
analogue of) \cite[Corollary 2.8.42]{Dal}.
It is possible that the converse implication holds as well. 
\par
Concluding this section, we shall see that $A_{\cb}(G)$ is
operator biprojective at least for those groups $G$ that satisfy the
hypotheses of Theorem \ref{thm1}.
\par
The key is the following lemma, whose straightforward proof we omit:
\begin{lemma} \label{tlem}
Let $E_1, E_2, F_1, F_2$ be operator spaces, and let $T_j \in \CB(E_j,F_j)$
be norm limits of finite rank operators for $j =1,2$. Then
$T_1 \tensor T_2 \in \CB(E_1 \Tensor E_2, F_1 \Tensor F_2)$ is a norm limit
of finite rank operators and thus compact.
\end{lemma}
\begin{proposition} \label{biprop}
Let $\A$ be a commutative, semisimple, Tauberian quantized Banach algebra
with discrete character space and 
a bounded approximate identity. Then $\A$ is operator biprojective
if and only if $\A$ is operator amenable. 
\end{proposition}
\begin{proof}
Any operator biprojective quantized Banach algebra with a bounded approximate
identity is operator amenable. Hence, only the ``if'' part needs proof.
\par
Suppose that $\A$ is operator amenable. By \cite[Proposition 2.4]{Rua}, 
it has an \emph{approximate diagonal}, i.e.\ a bounded net 
$( \boldsymbol{m}_\alpha )_{\alpha \in \mathbb A}$ such that
\[
  a \cdot \boldsymbol{m}_\alpha -\boldsymbol{m}_\alpha \cdot a \to 0
  \qquad (a \in \A)
\]
and 
\[
  a \Delta \boldsymbol{m}_\alpha \to a \qquad (a \in \A).
\]
For $a \in \A$, let $L_a, R_a \in \CB(\A)$ denote the operator of left and
right multiplication by $a$, respectively. Since $\A$ is semisimple and
Tauberian and has a discrete character space, $L_a$ and $R_a$ are norm
limits of finite rank operators for each $a \in \A$. 
Let $\cal U$ be an ultrafilter on $\mathbb A$
dominating the order filter, and let $a \in \A$. By Cohen's factorization
theorem (\cite[Theorem 2.9.24]{Dal}), 
there are $b,c \in \A$ such that $a = bc$, and by Lemma \ref{tlem},
$L_b \tensor R_c \in \CB(\A \Tensor \A)$ is compact. It follows that
\[
  \lim_{\cal U} a \cdot \boldsymbol{m}_\alpha =
  \lim_{\cal U} b \cdot \boldsymbol{m}_\alpha \cdot c =
  \lim_{\cal U} (L_b \tensor R_c) \boldsymbol{m}_\alpha
\]
exists. Define
\[
  \rho \!: \A \to \A \Tensor \A, \quad a \mapsto  
  \lim_{\cal U} a \cdot \boldsymbol{m}_\alpha.
\]
Then $\rho$ is completely bounded, and easily seen to be an $\A$-bimodule 
homomorphism and a right inverse of $\Delta$.
\end{proof}
\begin{remark}
The proof of the non-obvious direction of Proposition \ref{biprop}
is very similar to that of \cite[Corollary 3.2]{LRRW}. However, we do not
know if a straightforward quantization of \cite[Corollary 3.2]{LRRW} is
possible: unlike for the projective tensor product of Banach spaces, we
do not know whether the tensor product of two compact, completely bounded
maps between operator spaces is a compact map between the corresponding
projective tensor products (of operator spaces).
\end{remark}
\par
In view of Proposition \ref{prop1} and Theorem \ref{thm1}, we obtain:
\begin{corollary}
Let $G$ be a weakly amenable, discrete group such that $\cstar(G)$ is 
residually finite-dimensional. Then $A_{\cb}(G)$ is operator 
biprojective.
\end{corollary}
\section{Complementation of ideals in $A(G)$ and $A_{\cb}(G)$: an application}
In this section, we will consider (complete) complementation properties
of ideals in $A(G)$ and $A_{\cb}(G)$, where $G$ is a discrete group such that 
$A_{\cb}(G)$ is operator amenable.
\par
Let $G$ be a locally compact (mostly discrete) group, and let $F \subset G$
be closed. We set
\[
  I(F) := \{ f \in A(G) : f |_F \equiv 0 \}
  \qquad\text{and}\qquad
  I_{\cb}(E) := \{ f \in A_{\cb}(G) :  f |_F \equiv 0 \}.
\]
(Since $A(G)$ and $A_{\cb}(G)$ have the same character space, we use
different symbols when dealing with $A(G)$ and $A_{\cb}(G)$, 
respectively, in order to avoid confusion.) Similarly, we define 
\[
  J(F) := \varcl{\{ f \in A(G) : 
  \text{$\supp(f)$ is compact and has empty intersection with $F$} \}}^{\|
  \cdot \|_{A(G)}}
\]
and
\[
  J_{\cb}(F) := \varcl{\{ f \in A(G) : 
  \text{$\supp(f)$ is compact and has empty intersection with $F$} \}}^{\|
  \cdot \|_{A_{\cb}(G)}}.
\]
We say that $F$ is a \emph{set of synthesis} 
for $A(G)$ or $A_{\cb}(G)$, respectively, if $J(F)=I(F)$ or 
$J_{\cb}(F) = I_{\cb}(F)$, respectively.
\par
We begin with a useful observation:
\begin{proposition} \label{prop2}
Let $G$ be a weakly amenable locally compact group, and let $F \subset G$ be 
be a set of synthesis for $A(G)$. Then $F$ is a set of 
synthesis for $A_{\cb}(G)$.
\end{proposition} 
\begin{proof}
We first claim that $I(F)$ is dense in $I_{\cb}(F)$. To see this, let
$( e_\alpha )_{\alpha \in \mathbb A}$ be a bounded approximate identity for 
$A_{\cb}(G)$ contained in $A(G)$, and let $f \in I_{\cb}(G)$,
so that $f =\lim_\alpha f e_\alpha$. Since $f e_\alpha \in 
A(G)\cap I_{\cb}(F)= I(F)$ for each $\alpha \in \mathbb A$, this 
proves the claim.
\par
Since $F$ is a set of spectral synthesis for $A(G)$, it follows that
$I(F) \subset J_{\cb}(F)$, so that $I_{\cb}(F) =
J_{\cb}(F)$. 
\end{proof}
\begin{corollary} \label{cor2}
Let $G$ be a discrete and weakly amenable group. Then every subset of
$G$ is a set of synthesis for $A_{\cb}(G)$.
\end{corollary}
\begin{proof}
This follows immediately from the Proposition \ref{prop2} 
and (\cite[Proposition 2.2]{KL}).
\end{proof}
\par
Let $E$ be an operator space, and let $F$ be a closed subspace of $E$. 
We say that $F$ is \emph{completely complemented} in $E$ if there exists a 
completely bounded projection $P$ from $E$ onto $F$, and we say that
$F$ is \emph{completely weakly complemented} in $E$ if there exists a 
completely bounded projection from $E^\ast$ onto $F^\perp$. As in the
classical situation (\cite[Theorem 2.3.7]{LoA}), 
a closed ideal in an operator amenable, 
quantized Banach algebra is operator amenable if and only if it is weakly
complemented and if and only if it has a bounded approximate identity
(see \cite[Lemma 1.6]{RS}).
\par
The following proposition adds two more equivalent statements in the
case where the quantized Banach algebra is of the form $A_{\cb}(G)$
for a discrete group $G$:
\begin{proposition} \label{prop3}
Let $G$ be a discrete group such that $A_{\cb}(G)$ is operator
 amenable, and let $I$ be a closed ideal of $A_{\cb}(G)$.
Then the following are equivalent:
\begin{items}
\item $I$ is completely complemented;
\item $I$ is completely weakly complemented;
\item there is $F \subset G$ with $1_{F} \in  {\mathcal M}_{\cb}(A(G))$ such
that $I = I_{\cb}(F)$;
\item $I$ has an approximate identity bounded in the $\cb$-multiplier norm;
\item $I$ is operator amenable. 
\end{items}
\end{proposition}
\begin{proof}
As already stated, (ii) $\Longleftrightarrow$ (iv) $\Longleftrightarrow$ (v)
are well known (and hold for any closed ideal in a quantized Banach algebra).
Furthermore, (i) $\Longrightarrow$ (ii) is trivial.
\par
(iv) $\Longrightarrow$ (iii): Let $F$ be the hull of $I$, i.e.\
$F := \{ x \in G : \text{$f(x) = 0$ for all $f \in I$} \}$. By Corollary
\ref{cor2}, we have $I = I_{\cb}(F)$. Let 
$( e_\alpha )_\alpha$ be a bounded approximate identity for $I$. Since
${\cal M}_{\cb}(A(G))$ is a dual space, we can suppose that 
$( e_\alpha )_\alpha$ converges in the $w^\ast$-topology to some $f \in
{\cal M}_{\cb}(A(G))$. Since $w^\ast$-convergence in ${\cal M}_{\cb}(A(G))$
entails pointwise convergence on $G$, it follows that $f = 1_{G \setminus F}$,
so that $1_F = 1 - 1_{G \setminus F} \in {\cal M}_{\cb}(G)$.
\par
(iii) $\Longrightarrow$ (i): Since $1_F \in {\cal M}_{\cb}(A(G))$, the map
\[
  A_{\cb}(G) \to A_{\cb}(G), \quad f \mapsto 1_{G \setminus F} f
\]
is a completely bounded projection onto $I$.
\end{proof}
\begin{remark}
The first four equivalences of Proposition \ref{prop3} can be viewed as 
extensions of the main results of (\cite{Woo1}), which were primarily about
the Fourier algebra of an amenable group. The proof of Proposition 
\ref{prop3} is very similar to the corresponding arguments in 
\cite{Woo1}.
\end{remark}
\par
A somewhat more surprising result is that, under the same hypotheses
as in Proposition \ref{prop3}, we can obtain the equivalence of
(i) to (iv) for closed ideals of the Fourier algebra with its original norm
(Corollary \ref{cor3}, below).
\par
The crucial implication is the following:
\begin{theorem} \label{thm2}
Let $G$ be a discrete group such that $A_{\cb}(G)$ is operator
amenable, and let $I$ be a weakly completely complemented closed ideal of 
$A(G)$. Then there is $F \subset G$ with $1_F \in {\cal M}_{\cb}(G)$ such that
$I=I(F)$.
\end{theorem}
\begin{proof}
Since $G$ is discrete and weakly amenable, \cite[Proposition 2.2]{KL} yields
$F \subset G$ such that $I=I(F)$. It remains to be shown that $1_F
\in {\cal M}_{\cb}(G)$.
\par
Since $I$ is an ideal of $A(G)$, it is a weakly completely complemented
$A_{\cb}(G)$-submodule of the (symmetric) quantized Banach 
$A_{\cb}(G)$-module $A(G)$. By definition, $I^\perp$ thus 
is a completely complemented, closed $A_{\cb}(G)$-submodule of the
dual $A_{\cb}(G)$-module $\VN(G)$. Since $A_{\cb}(G)$ is 
operator amenable, \cite[Theorem 1]{Woo1} implies that $I^{\perp }$ is 
completely invariantly complemented, i.e.\ there is a completely bounded 
projection $P \!: \VN(G)\to I^\perp$ which is an $A_{\cb}(G)$-module
homomorphism.
\par
Define $Q \!: A(G)^{\ast\ast} \to A(G)^{\ast\ast}$ as the complementary
projection of $P^\ast$, i.e.\ $Q := \id_{A(G)^{\ast\ast}} -P^\ast$.
Then $Q$ is a completely bounded projection from 
$A(G)^{\ast\ast}$ onto $(I^\perp)^\perp = I^{\ast\ast}$ and an 
$A_{\cb}(G)$-module homomorphism.
\par
Let $x \in G$, so that $1_{\{ x \}} \in A(G)$, and note that
\[
  Q\left(1_{ \{ x \}}\right) = 
  Q\left(1_{ \{ x \}}^2\right) = 
  1_{ \{ x \}} \cdot Q\left(1_{ \{ x \}}\right).
\]
Since $G$ is discrete, $A(G)$ is an ideal in $A(G)^{\ast\ast}$, so that
$Q(1_{ \{ x \}}) \in A(G)$. Since $A(G)$ is Tauberian, it follows that
$Q(A(G)) \subset A(G)$.
\par
All in all, $Q$ is completely bounded, maps $A(G)$ into itself, and
is an $A_{\cb}(G)$-module homomorphism. It follows that $Q |_{A(G)}$ is 
a completely bounded multiplier of $A(G)$, i.e.\ there is $g \in 
{\mathcal M}_{\cb}(A(G))$ such that $Q f = gf$ for all $f \in A(G)$. Finally,
as $Q$ is a projection onto $I^{\ast\ast}$, it is clear that 
$g = 1_{G \setminus F}$, so that $1_F \in  {\mathcal M}_{\cb}(A(G))$.
\end{proof}
\begin{corollary} \label{cor3}
Let $G$ be a discrete group such that $A_{\cb}(G)$ is operator
amenable, and let $I$ be a closed ideal of $A(G)$. Then the following are 
equivalent:
\begin{items}
\item $I$ is completely complemented;
\item $I$ is completely weakly complemented;
\item there is $F \subset G$ with $1_{F} \in  {\mathcal M}_{\cb}(A(G))$ such
that $I = I(F)$;
\item $I$ has an approximate identity bounded in the $\cb$-multiplier norm.
\end{items}
\end{corollary}
\begin{proof}
(i) $\Longrightarrow$ (ii) is trivial, (ii) $\Longrightarrow$ (iii)
follows from Theorem \ref{thm2}, and (iii) $\Longrightarrow$ (i) follows
as in the proof of Proposition \ref{prop3}.
\par
(iii) $\Longrightarrow$ (iv): 
Since $A_{\cb}(G)$ is operator amenable, it has a bounded approximate 
identity, so that $A(G)$ has an approximate identity, say
 $( e_\alpha )_\alpha$, that is bounded in $A_{\cb}(G)$
(\cite[Proposition 1]{For}). Then $( 1_{G \setminus F} e_\alpha )_\alpha$
is the desired approximate identity.
\par
(iv) $\Longrightarrow$ (iii): This is proven as the corresponding implication
of Proposition \ref{prop3}.
\end{proof}
\begin{remark}
The equivalence of Corollary \ref{cor3}(i) and (iii) was proven by 
Wood, first for amenable discrete groups in \cite{Woo1} and then, later,
for all discrete groups in (\cite{Woo2}). Wood's techniques, however, do not
allow to prove the equivalence of (i) and (ii) or of (i) and (iii) with (iv) 
without the stronger hypothesis that $G$ be amenable.
\end{remark}
\renewcommand{\baselinestretch}{1.0}
\vfill
\begin{tabbing}
{\it Second author's address\/}: \= Department of Mathematical and Statistical Sciences \kill 
{\it First author's address\/}:  \> Department of Pure Mathematics \\
                                  \> University of Waterloo \\
                                  \> Waterloo, Ontario \\
                                  \> Canada N2L 3G1 \\[\medskipamount]
{\it E-mail\/}:                   \> {\tt beforres@math.ualberta.ca}\\[\bigskipamount]
{\it Second author's address\/}: \> Department of Mathematical and Statistical Sciences \\
                                  \> University of Alberta \\
                                  \> Edmonton, Alberta \\
                                  \> Canada T6G 2G1 \\[\medskipamount]
{\it E-mail\/}:                   \> {\tt vrunde@ualberta.ca}\\[\medskipamount]
{\it URL\/}:                      \> {\tt http://www.math.ualberta.ca/$^\sim$runde/} \\[\bigskipamount]
{\it Third author's address\/}:  \> Department of Pure Mathematics \\
                                  \> University of Waterloo \\
                                  \> Waterloo, Ontario \\
                                  \> Canada N2L 3G1 \\[\medskipamount]
{\it E-mail\/}:                   \> {\tt nspronk@math.ualberta.ca}\\[\medskipamount]
{\it URL\/}:                      \> {\tt http://www.math.uwaterloo.ca/$^\sim$nspronk/} 
\end{tabbing}
\end{document}